\documentclass[11pt]{article}
\usepackage{epsfig}

\def\be{\begin{equation}}
\def\ee{\end{equation}}
\def\bea{\begin{eqnarray}}
\def\eea{\end{eqnarray}}
\def\bes{\begin{eqnarray*}}
\def\ees{\end{eqnarray*}}

\def\nn{\nonumber}
\def\<{\langle}
\def\>{\rangle}
\def\lb{\label}
\def\bs{\setminus}

\def\R{{\bf R}}
\def\C{{\bf C}}
\def\Z{{\bf Z}}
\def\N{{\bf N}}
\def\U{{\bf U}}
\def\Q{{\bf Q}}

\def\ga{{\gamma}}

\def\th{{\theta}}
\def\om{{\omega}}
\def\Om{{\Omega}}
\def\ep{{\epsilon}}
\def\lm{{\lambda}}

\def\sg{{\sigma}}

\def\Sg{{\Sigma}}

\def\H{{\cal H}}

\def\P{{\cal P}}
\def\J{{\cal J}}

\def\Nn{{\cal N}}

\def\im{{\rm im}}

\def\Sp{{\rm Sp}}

\def\dm{{\rm \diamond}}

\def\ol#1{\overline{#1}}  
\def\td#1{\tilde{#1}}
\def\hb{\vrule height0.18cm width0.14cm $\,$}

\def\ol#1{\overline{#1}}  
\def\td#1{\tilde{#1}}

\title{Stability of closed characteristics on symmetric \\compact convex hypersurfaces in $\R^{2n}$}

\author{Wei Wang\thanks{Partially supported by LMAM in Peking University
in China and China Postdoctoral Science Foundation No.20070420264.
E-mail: alexanderweiwang@yahoo.com.cn, wangwei@math.pku.edu.cn  }\\
School of Mathematical Science \\ Peking University, Beijing 100871 \\
PEOPLES REPUBLIC OF CHINA \\ }
\date{}

\date{}

\begin{document}

\maketitle

\begin{abstract}
{\it In this article, let $\Sigma\subset\R^{2n}$ be a compact convex
hypersurface which is symmetric with respect to the origin. We prove
that if $\Sg$ carries finitely many geometrically distinct closed
characteristics, then at least $n-1$ of them must be non-hyperbolic;
if $\Sg$ carries exactly $n$ geometrically distinct closed
characteristics, then at least two of them must be elliptic. }
\end{abstract}



{\bf Running title}: Stability of closed characteristics

\renewcommand{\theequation}{\thesection.\arabic{equation}}
\renewcommand{\thefigure}{\thesection.\arabic{figure}}

\setcounter{equation}{0}
\section{Introduction and main results}

In this article, let $\Sigma$ be a fixed $C^3$ compact convex hypersurface
in $\R^{2n}$, i.e., $\Sigma$ is the boundary of a compact and strictly
convex region $U$ in $\R^{2n}$. We denote the set of all such hypersurfaces
by $\H(2n)$. Without loss of generality, we suppose $U$ contains the origin.
We denote the set of all compact convex hypersurfaces which are
symmetric with respect to the origin by $\mathcal{SH}(2n)$,
i.e., $\Sg=-\Sg$ for $\Sg\in\mathcal{SH}(2n)$.
We consider closed characteristics $(\tau,y)$ on $\Sigma$, which are
solutions of the following problem
\be \left\{\matrix{\dot{y}=JN_{\Sigma}(y), \cr
               y(\tau)=y(0), \cr }\right. \lb{1.1}\ee
where $J=\left(\matrix{0 &-I_n\cr
                I_n  & 0\cr}\right)$,
$I_n$ is the identity matrix in $\R^n$, $\tau>0$ and $N_\Sigma(y)$ is
the outward normal vector of $\Sigma$ at $y$ normalized by the
condition $N_{\Sigma}(y)\cdot y=1$. Here $a\cdot b$ denotes the
standard inner product of $a, b\in\R^{2n}$. A closed characteristic
$(\tau,\, y)$ is {\it prime} if $\tau$ is the minimal period of $y$.
Two closed characteristics $(\tau,\, y)$ and $(\sigma, z)$ are {\it
geometrically distinct}  if $y(\R)\not= z(\R)$. We denote by
$\J(\Sg)$ and $\widetilde{\J}(\Sg)$ the set of all closed
characteristics $(\tau,\, y)$ on $\Sg$ with $\tau$ being the minimal
period of $y$ and the set of all geometrically distinct ones
respectively. Note that
$\J(\Sg)=\{\theta\cdot y\,|\, \theta\in S^1,\;y\; is\; prime\}$,
while $\widetilde{\J}(\Sg)=\J(\Sg)/S^1$, where the natural $S^1$-action
is defined by $\theta\cdot y(t)=y(t+\tau\theta),\;\;\forall \theta\in S^1,\,t\in\R$.

Let $j: \R^{2n}\rightarrow\R$ be the gauge function of $\Sigma$, i.e.,
$j(\lambda x)=\lambda$ for $x\in\Sigma$ and $\lambda\ge0$, then
$j\in C^3(\R^{2n}\setminus\{0\}, \R)\cap C^0(\R^{2n}, \R)$
and $\Sigma=j^{-1}(1)$. Fix a constant $\alpha\in(1,\,2)$ and
define the Hamiltonian function
$H_\alpha :\R^{2n}\rightarrow [0,\,+\infty)$ by
\be H_\alpha(x)=j(x)^\alpha,\qquad \forall x\in\R^{2n}.\lb{1.2}\ee
Then
$H_\alpha\in C^3(\R^{2n}\setminus\{0\}, \R)\cap C^1(\R^{2n}, \R)$
is convex and $\Sigma=H_\alpha^{-1}(1)$.
It is well known that problem (\ref{1.1}) is equivalent to
the following given energy problem of Hamiltonian system
\be
\left\{\matrix{\dot{y}(t)=JH_\alpha^\prime(y(t)),
             &&\quad H_\alpha(y(t))=1,\qquad \forall t\in\R. \cr
     y(\tau)=y(0). && \cr }\right. \lb{1.3}\ee
Denote by $\mathcal{J}(\Sigma, \,\alpha)$ the set of all solutions
$(\tau,\, y)$ of (\ref{1.3}) where $\tau$ is the minimal period of
$y$ and by $\widetilde{\mathcal{J}}(\Sigma, \,\alpha)$ the set of
all geometrically distinct solutions  of (\ref{1.3}). As above,
$\widetilde{\mathcal{J}}(\Sigma, \,\alpha)$ is obtained from
$\mathcal{J}(\Sigma, \,\alpha)$ by dividing the natural
$S^1$-action. Note that elements in $\mathcal{J}(\Sigma)$ and
$\mathcal{J}(\Sigma, \,\alpha)$ are one to one correspondent to each
other, similarly for $\widetilde{\J}(\Sg)$ and
$\widetilde{\mathcal{J}}(\Sigma, \,\alpha)$.

Let $(\tau,\, y)\in\mathcal{J}(\Sigma, \,\alpha)$. The fundamental
solution $\gamma_y : [0,\,\tau]\rightarrow \Sp(2n)$ with $\gamma_y(0)=I_{2n}$
of the linearized Hamiltonian system
\be \dot w(t)=JH_\alpha^{\prime\prime}(y(t))w(t),\qquad \forall t\in\R,\lb{1.4}\ee
is called the {\it associate symplectic path} of $(\tau,\, y)$.
The eigenvalues of $\gamma_y(\tau)$ are called {\it Floquet multipliers}
of $(\tau,\, y)$. By Proposition 1.6.13 of \cite{Eke3}, the Floquet multipliers
with their multiplicities of $(\tau,\, y)\in\mathcal{J}(\Sigma)$ do not depend on
the particular choice of the Hamiltonian function in (\ref{1.3}).
For any $M\in \Sp(2n)$, we define the {\it elliptic height } $e(M)$ of
$M$ to be the total algebraic multiplicity of all eigenvalues of $M$ on the
unit circle $\U=\{z\in\C|\; |z|=1\}$ in the complex plane $\C$.
Since $M$ is symplectic, $e(M)$ is even and $0\le e(M)\le 2n$.
As usual $(\tau,\, y)\in\J(\Sg)$ is {\it elliptic} if
$e(\gamma_y(\tau))=2n$. It is {\it non-degenerate} if $1$ is a double
Floquet multiplier of it. It is  {\it hyperbolic} if $1$ is a
double Floquet multiplier of it and $e(\gamma_y(\tau))=2$.
It is {\it irrationally elliptic} if $\gamma_y(\tau)$ is suitably
homotopic to the $\dm$-product of one $N_1(1,1)=\left(\matrix{1 & 1 \cr
                   0 & 1 \cr}\right)$ and $n-1$ rotation $2\times 2$ matrices
with rotation angles being irrational multiples of $\pi$, more precisely,
$N_1(1,1)\dm R(\th_1)\dm\cdots\dm R(\th_{n-1})\in \Om^0(\ga_x(\tau))$ for some
$\th_i\in (0,2\pi)\bs (\pi\Q)$ with $1\le i\le n-1$, cf., \S3 below for notations.
It is well known that these concepts are independent of the choice of $\alpha\in(1,\,2)$.

For the existence and multiplicity of geometrically distinct closed
characteristics on convex compact hypersurfaces in $\R^{2n}$ we refer to
\cite{Rab1}, \cite{Wei1}, \cite{EkL1}, \cite{EkH1}, \cite{Szu1},
\cite{Vit1}, \cite{HWZ}, \cite{LoZ1}, \cite{LLZ}, \cite{WHL}, and references therein.

On the stability problem, in \cite{Eke2} of Ekeland in 1986 and \cite{Lon2}
of Long in 1998, for any $\Sg\in\H(2n)$ the existence of at least one
non-hyperbolic closed characteristic on $\Sg$ was proved provided
$^\#\td{\J}(\Sg)<+\infty$. Ekeland proved also in \cite{Eke2} the existence
of at least one elliptic closed characteristic on $\Sg$ provided $\Sg\in\H(2n)$
is $\sqrt{2}$-pinched. In \cite{DDE1} of 1992, Dell'Antonio, D'Onofrio and
Ekeland proved the existence of at least one elliptic closed characteristic
on $\Sg$ provided $\Sg\in\mathcal{SH}(2n)$. In \cite{Lon4} of 2000,
Long proved that $\Sg\in\mathcal{\H}(4)$ and $\,^{\#}\td{\J}(\Sg)=2$ imply that both of
the closed characteristics must be elliptic. In \cite{LoZ1} of 2002, Long and
Zhu further proved when $^\#\td{\J}(\Sg)<+\infty$, there exists at least one
elliptic closed characteristic and there are at least $[\frac{n}{2}]$ geometrically
distinct closed characteristics on $\Sg$ possessing irrational mean indices,
which are then non-hyperbolic. In \cite{LoW1}, Long and the author
proved that there exist at least two non-hyperbolic closed characteristics on
$\Sg\in\H(6)$ when $^\#\td{\J}(\Sg)<+\infty$.
In \cite{Wang}, the author proved that on every $\Sg\in\H(6)$ satisfying
$^\#\td{\J}(\Sg)<+\infty$, there exist at least two
closed characteristics possessing irrational mean indices and if
$^\#\td{\J}(\Sigma)=3$, then there exist at least two elliptic closed
characteristics. It was conjectured by Hofer et al. that
$ \{\,^{\#}\td{\J}(\Sg)\,|\,\Sg\in {\cal H}(2n)\} = \{n\}\cup\{+\infty\}$
for $n\ge 2$ and it was conjectured by Long et al. in \cite{WHL} that
all the closed characteristics
on $\Sg$ are irrationally elliptic for $\Sg\in {\cal H}(2n)$ with $n\ge 2$
whenever $\,^{\#}\td{\J}(\Sg)<\infty$.
Note that both conjectures have been proved in the $n=2$ case,
cf. \cite{HWZ} and \cite{WHL} respectively.

Motivated by these results, we prove
the following results in this article:

{\bf Theorem 1.1.} {\it On every $\Sg\in\mathcal{SH}(2n)$ satisfying
$^\#\td{\J}(\Sg)<+\infty$ there exist at least $n-1$
non-hyperbolic closed characteristics in $\td{\J}(\Sigma)$. }

{\bf Theorem 1.2.} {\it Suppose $^\#\td{\J}(\Sigma)=n$ for some
$\Sigma\in\mathcal{SH}(2n)$ and $n\ge 2$. Then there exist at least two elliptic closed
characteristics in $\td{\J}(\Sigma)$.}

The proofs of these theorems are given in \S4.
The proofs are motivated by the methods in \cite{LoZ1}
and \cite{LLZ} by using the index iteration theory developed
by Long and his coworkers, specially the common index jump theorem of Long and Zhu
(Theorem 4.3 of \cite{LoZ1}, cf. Theorem 11.2.1 of \cite{Lon5}).
In \S2 and \S3, we review briefly the variational structure
for closed characteristics and the index iteration theory for
symplectic paths respectively.

In this article, let $\N$, $\N_0$, $\Z$, $\Q$, $\R$, and $\C$ denote
the sets of natural integers, non-negative integers, integers, rational
numbers, real numbers, and complex numbers respectively.
Denote by $a\cdot b$ and $|a|$ the standard inner product and norm in
$\R^{2n}$. Denote by $\langle\cdot,\cdot\rangle$ and $\|\cdot\|$
the standard $L^2$-inner product and $L^2$-norm. For an $S^1$-space $X$, we denote
by $X_{S^1}$ the homotopy quotient of $X$ module the $S^1$-action, i.e.,
$X_{S^1}=S^\infty\times_{S^1}X$. We define the functions
\be \left\{\matrix{[a]=\max\{k\in\Z\,|\,k\le a\}, &
E(a)=\min\{k\in\Z\,|\,k\ge a\} , \cr
                   \varphi(a)=E(a)-[a],   \cr}\right. \lb{1.5}\ee
Specially, $\varphi(a)=0$ if $ a\in\Z\,$, and $\varphi(a)=1$ if $a\notin\Z\,$.
In this article we use only $\Q$-coefficients for all homological modules.

\setcounter{equation}{0}
\section{Variational structure for closed characteristics }

In this section, we describe the variational structure for closed
characteristics.

As in P.199 of \cite{Eke3}, choose some $\alpha\in(1,\, 2)$ and associate with $U$
a convex function $H_\alpha$ such that $H_\alpha(\lambda x)=\lambda^\alpha H_\alpha(x)$ for $\lambda\ge 0$.
Consider the fixed period problem
\be \left\{\matrix{\dot{x}(t)=JH_\alpha^\prime(x(t)), \cr
     x(1)=x(0).         \cr }\right. \lb{2.1}\ee

Define
\be L_0^{\frac{\alpha}{\alpha-1}}(S^1,\R^{2n})
  =\{u\in L^{\frac{\alpha}{\alpha-1}}(S^1,\R^{2n})\,|\,\int_0^1udt=0\}. \lb{2.2}\ee
The corresponding Clarke-Ekeland dual action functional is defined by
\be \Phi(u)=\int_0^1\left(\frac{1}{2}Ju\cdot Mu+H_\alpha^{\ast}(-Ju)\right)dt,
    \qquad \forall\;u\in L_0^{\frac{\alpha}{\alpha-1}}(S^1,\R^{2n}), \lb{2.3}\ee
where $Mu$ is defined by $\frac{d}{dt}Mu(t)=u(t)$ and $\int_0^1Mu(t)dt=0$,
$H_\alpha^\ast$ is the Fenchel transform of $H_\alpha$ defined by
$H_\alpha^\ast(y)=\sup\{x\cdot y-H_\alpha(x)\;|\; x\in \R^{2n}\}$.
By Theorem 5.2.8 of \cite{Eke3}, $\Phi$ is $C^1$ on
$L_0^{\frac{\alpha}{\alpha-1}}$ and satisfies
the Palais-Smale condition.  Suppose
$x$ is a solution of (\ref{2.1}). Then $u=\dot{x}$ is a critical point
of $\Phi$. Conversely, suppose $u$ is a critical point of $\Phi$.
Then there exists a unique $\xi\in\R^{2n}$ such that $Mu-\xi$ is a
solution of (\ref{2.1}). In particular, solutions of (\ref{2.1}) are in
one to one correspondence with critical points of $\Phi$. Moreover,
$\Phi(u)<0$ for every critical point $u\not= 0$ of $\Phi$.

Suppose $u$ is a nonzero critical point of $\Phi$. Then
the formal Hessian of $\Phi$ at $u$ on $L_0^2(S^1,\R^{2n})$ is defined by
$$ Q(v,\; v)=\int_0^1 (Jv\cdot Mv+(H_\alpha^\ast)^{\prime\prime}(-Ju)Jv\cdot Jv)dt, $$
which defines an orthogonal splitting $L_0^2(S^1,\R^{2n})=E_-\oplus E_0\oplus E_+$ of
$L_0^2(S^1,\; \R^{2n})$ into negative, zero and positive subspaces. The
index of $u$ is defined by $i(u)=\dim E_-$ and the nullity of $u$ is
defined by $\nu(u)=\dim E_0$. Specially $1\le \nu(u)\le 2n$
always holds, cf. P.219 of  \cite{Eke3}.

We have a natural $S^1$-action on $L_0^{\frac{\alpha}{\alpha-1}}(S^1,\; \R^{2n})$ defined by
$\th\cdot u(t)=u(\th+t)$ for all $\th\in S^1$ and $t\in\R$. Clearly
$\Phi$ is $S^1$-invariant. Hence if $u$ is a critical
point of $\Phi$, then the whole orbit $S^1\cdot u$ is formed by
critical points of $\Phi$. Denote by $crit(\Phi)$ the set of
critical points of $\Phi$.
Then we make the following definition

{\bf Definition 2.1.} {\it Suppose $u$ is a nonzero critical
point of $\Phi$, and $\Nn$ is an $S^1$-invariant
open neighborhood of $S^1\cdot u$ such that
$crit(\Phi)\cap(\Lambda(u)\cap \Nn)=S^1\cdot u$. Then
the $S^1$-critical modules of $S^1\cdot u$ is defined by
\bea C_{S^1,\; q}(\Phi, \;S^1\cdot u)
=H_{q}((\Lambda(u)\cap\Nn)_{S^1},\;
((\Lambda(u)\setminus S^1\cdot u)\cap\Nn)_{S^1}),\lb{2.4}
\eea
where $\Lambda(u)=\{w\in L_0^{\frac{\alpha}{\alpha-1}}(S^1,\R^{2n})\;|\;
\Phi(w)\le\Phi(u)\}$.}

By the proof of Proposition 3.6 of \cite{Wang} we have
\be C_{S^1,\; \ast}(\Phi, \;S^1\cdot u)\cong C_{S^1,\; \ast}(\Psi_a, \;S^1\cdot u_a),\lb{2.5}\ee
where $\Psi_a$ is the functional constructed in \cite{WHL} and $u_a$ is its
critical point corresponding to $u$.
By Proposition 3.5 of \cite{WHL}, the index and nullity of $\Psi_a$
at $u_a$ coincide with those of $\Phi$ at $u$.
Hence by Propositions 2.3 and 2.6 of \cite{Wang}, we have

{\bf Proposition 2.2.} {\it Let $k_j(u)\equiv\dim C_{S^1,\; j}(\Phi, \;S^1\cdot u)$.
Then $k_j(u)$  equal to $0$ when $j<i(u)$ or $j>i(u)+\nu(u)-1$ and can only take values $0$ or $1$
when $j=i(u)$ or $j=i(u)+\nu(u)-1$. }

For a closed characteristic $(\tau,y)$ on $\Sigma$, we denote by
$y^m\equiv (m\tau, y)$ the $m$-th iteration of $y$ for $m\in\N$.
Let $u^m$ be the unique critical point of $\Phi$ corresponding to $(m\tau, y)$.
Then we define the index $i(y^m)$ and nullity $\nu(y^m)$
of $(m\tau,y)$ for $m\in\N$ by
\be i(y^m)=i(u^m), \qquad \nu(y^m)=\nu(u^m). \lb{2.6}\ee
The mean index of $(\tau,y)$ is defined by
\be \hat{i}(y)=\lim_{m\rightarrow\infty}\frac{i(y^m)}{m}. \lb{2.7}\ee
Note that $\hat{i}(y)>2$ always holds which was proved by Ekeland and
Hofer in \cite{EkH1} of 1987 (cf. Corollary 8.3.2 and Lemma 15.3.2
of \cite{Lon5} for a different proof).

Recall that for a principal $U(1)$-bundle $E\to B$, the Fadell-Rabinowitz index
(cf. \cite{FaR1}) of $E$ is defined to be $\sup\{k\;|\, c_1(E)^{k-1}\not= 0\}$,
where $c_1(E)\in H^2(B,\Q)$ is the first rational Chern class. For a $U(1)$-space,
i.e., a topological space $X$ with a $U(1)$-action, the Fadell-Rabinowitz index is
defined to be the index of the bundle $X\times S^{\infty}\to X\times_{U(1)}S^{\infty}$,
where $S^{\infty}\to CP^{\infty}$ is the universal $U(1)$-bundle.

For any $\kappa\in\R$, we denote by
\be \Phi^{\kappa-}=\{u\in L_0^{\frac{\alpha}{\alpha-1}}(S^1,\R^{2n})\;|\;
             \Phi(u)<\kappa\}. \lb{2.8}\ee
Then as in P.218 of \cite{Eke3}, we define
\be c_i=\inf\{\delta\in\R\;|\: \hat I(\Phi^{\delta-})\ge i\},\lb{2.9}\ee
where $\hat I$ is the Fadell-Rabinowitz index defined above. Then by Proposition 3
in P.218 of \cite{Eke3}, we have

{\bf Proposition 2.3.} {\it Every $c_i$ is a critical value of $\Phi$. If
$c_i=c_j$ for some $i<j$, then there are infinitely many geometrically
distinct closed characteristics on $\Sg$.}

Comparing with Theorem 4 in P.219 of \cite{Eke3}, we have the following
property by Proposition 3.5 of \cite{Wang}.

{\bf Proposition 2.4.} {\it For every $i\in\N$, there exists a point
$u\in L_0^{\frac{\alpha}{\alpha-1}}(S^1,\R^{2n})$ such that}
\bea
&& \Phi^\prime(u)=0,\quad \Phi(u)=c_i, \lb{2.10}\\
&& C_{S^1,\; 2(i-1)}(\Phi, \;S^1\cdot u)\neq 0. \lb{2.11}\eea

{\bf Definition 2.5.} {\it A prime closed characteristic $(\tau,\, y)$ is
$(m, i)$-{ \bf variationally visible:} if there exist some $m, i\in\N$
such that (\ref{2.10}) and (\ref{2.11}) hold for $y^m$ and $c_i$. We call $(\tau,\, y)$
{\bf infinitely variationally visible:} if there exist
infinitely many $m, i\in\N$ such that $(\tau,\, y)$ is $(m, i)$-variationally visible.
We denote by $\mathcal{V}(\Sg, \alpha)$ and $\mathcal{V}_\infty(\Sg, \alpha)$ the set of
variationally visible and infinitely
variationally visible closed characteristics respectively.}

Recall that the action of a closed characteristic $(\tau,\, y)$
is defined by (cf. P190 of \cite{Eke3})
\be A(\tau, y)=\frac{1}{2}\int_0^\tau (Jy\cdot \dot y)dt. \lb{2.12}\ee

Then we have the following

{\bf Theorem 2.6.} {\it Suppose there are only finitely many prime closed
 characteristics on $\Sg$. Then for any $(\tau, y)\in\mathcal{V}_\infty(\Sg, \alpha)$,
we have
\be \frac{\hat i(y)}{A(\tau, y)}=\frac{1}{\gamma(\Sg)}, \lb{2.13}
\ee
where
$$\gamma(\Sg)=C_\alpha^{-1}\lim\inf_{i\rightarrow\infty}(i(-c_i)^{\frac{2-\alpha}{\alpha}})^{-1}=
C_\alpha^{-1}\lim\sup_{i\rightarrow\infty}(i(-c_i)^{\frac{2-\alpha}{\alpha}})^{-1}$$
and $C_\alpha=
\frac{4}{\alpha}(1-\frac{\alpha}{2})^{\frac{\alpha-2}{\alpha}}$
}

{\bf Proof.} Note that we have ${\hat i(y^m)}=m{\hat i(y)}$ by
(\ref{2.7}) and $A(y^m)=mA(y)$ by (\ref{2.12}).
Thus $\frac{\hat i(y^m)}{A(y^m)}=\frac{\hat i(y)}{A(y)}$ for any
$m\in\N$. Now the theorem follows from Lemma 5.3.12 and
Theorem 5.3.15 of \cite{Eke3}.\hfill\hb

\setcounter{equation}{0}
\section{ A brief review on an index theory for symplectic pathss}

In this section, we recall briefly an index theory for symplectic paths
developed by Y. Long and his coworkers.
All the details can be found in \cite{Lon5}.

As usual, the symplectic group $\Sp(2n)$ is defined by
$$ \Sp(2n) = \{M\in {\rm GL}(2n,\R)\,|\,M^TJM=J\}, $$
whose topology is induced from that of $\R^{4n^2}$. For $\tau>0$ we are interested
in paths in $\Sp(2n)$:
$$ \P_{\tau}(2n) = \{\ga\in C([0,\tau],\Sp(2n))\,|\,\ga(0)=I_{2n}\}, $$
which is equipped with the topology induced from that of $\Sp(2n)$. The
following real function was introduced in \cite{Lon3}:
$$ D_{\om}(M) = (-1)^{n-1}\ol{\om}^n\det(M-\om I_{2n}), \qquad
          \forall \om\in\U,\, M\in\Sp(2n). $$
Thus for any $\om\in\U$ the following codimension $1$ hypersurface in $\Sp(2n)$ is
defined in \cite{Lon3}:
$$ \Sp(2n)_{\om}^0 = \{M\in\Sp(2n)\,|\, D_{\om}(M)=0\}.  $$
For any $M\in \Sp(2n)_{\om}^0$, we define a co-orientation of $\Sp(2n)_{\om}^0$
at $M$ by the positive direction $\frac{d}{dt}Me^{t\ep J}|_{t=0}$ of
the path $Me^{t\ep J}$ with $0\le t\le 1$ and $\ep>0$ being sufficiently
small. Let
\bea
\Sp(2n)_{\om}^{\ast} &=& \Sp(2n)\bs \Sp(2n)_{\om}^0,   \nn\\
\P_{\tau,\om}^{\ast}(2n) &=&
      \{\ga\in\P_{\tau}(2n)\,|\,\ga(\tau)\in\Sp(2n)_{\om}^{\ast}\}, \nn\\
\P_{\tau,\om}^0(2n) &=& \P_{\tau}(2n)\bs  \P_{\tau,\om}^{\ast}(2n).  \nn\eea
For any two continuous arcs $\xi$ and $\eta:[0,\tau]\to\Sp(2n)$ with
$\xi(\tau)=\eta(0)$, it is defined as usual:
$$ \eta\ast\xi(t) = \left\{\matrix{
            \xi(2t), & \quad {\rm if}\;0\le t\le \tau/2, \cr
            \eta(2t-\tau), & \quad {\rm if}\; \tau/2\le t\le \tau. \cr}\right. $$
Given any two $2m_k\times 2m_k$ matrices of square block form
$M_k=\left(\matrix{A_k&B_k\cr
                                C_k&D_k\cr}\right)$ with $k=1, 2$,
as in \cite{Lon5}, the $\;\dm$-product of $M_1$ and $M_2$ is defined by
the following $2(m_1+m_2)\times 2(m_1+m_2)$ matrix $M_1\dm M_2$:
$$ M_1\dm M_2=\left(\matrix{A_1&  0&B_1&  0\cr
                               0&A_2&  0&B_2\cr
                             C_1&  0&D_1&  0\cr
                               0&C_2&  0&D_2\cr}\right). \nn$$  
Denote by $M^{\dm k}$ the $k$-fold $\dm$-product $M\dm\cdots\dm M$. Note
that the $\dm$-product of any two symplectic matrices is symplectic. For any two
paths $\ga_j\in\P_{\tau}(2n_j)$ with $j=0$ and $1$, let
$\ga_0\dm\ga_1(t)= \ga_0(t)\dm\ga_1(t)$ for all $t\in [0,\tau]$.

A special path $\xi_n\in\P_{\tau}(2n)$ is defined by
\be \xi_n(t) = \left(\matrix{2-\frac{t}{\tau} & 0 \cr
                                             0 &  (2-\frac{t}{\tau})^{-1}\cr}\right)^{\dm n}
        \qquad {\rm for}\;0\le t\le \tau.  \lb{3.1}\ee

{\bf Definition 3.1.} (cf. \cite{Lon3}, \cite{Lon5}) {\it For any $\om\in\U$ and
$M\in \Sp(2n)$, define
\be  \nu_{\om}(M)=\dim_{\C}\ker_{\C}(M - \om I_{2n}).  \lb{3.2}\ee
For any $\tau>0$ and $\ga\in \P_{\tau}(2n)$, define
\be  \nu_{\om}(\ga)= \nu_{\om}(\ga(\tau)).  \lb{3.3}\ee

If $\ga\in\P_{\tau,\om}^{\ast}(2n)$, define
\be i_{\om}(\ga) = [\Sp(2n)_{\om}^0: \ga\ast\xi_n],  \lb{3.4}\ee
where the right hand side of (\ref{3.4}) is the usual homotopy intersection
number, and the orientation of $\ga\ast\xi_n$ is its positive time direction under
homotopy with fixed end points.

If $\ga\in\P_{\tau,\om}^0(2n)$, we let $\mathcal{F}(\ga)$
be the set of all open neighborhoods of $\ga$ in $\P_{\tau}(2n)$, and define
\be i_{\om}(\ga) = \sup_{U\in\mathcal{F}(\ga)}\inf\{i_{\om}(\beta)\,|\,
                       \beta\in U\cap\P_{\tau,\om}^{\ast}(2n)\}.
               \lb{3.5}\ee
Then
$$ (i_{\om}(\ga), \nu_{\om}(\ga)) \in \Z\times \{0,1,\ldots,2n\}, $$
is called the index function of $\ga$ at $\om$. }

Note that when $\om=1$, this index theory was introduced by
C. Conley-E. Zehnder in \cite{CoZ1} for the non-degenerate case with $n\ge 2$,
Y. Long-E. Zehnder in \cite{LZe1} for the non-degenerate case with $n=1$,
and Y. Long in \cite{Lon1} and C. Viterbo in \cite{Vit2} independently for
the degenerate case. The case for general $\om\in\U$ was defined by Y. Long
in \cite{Lon3} in order to study the index iteration theory (cf. \cite{Lon5}
for more details and references).

For any symplectic path $\ga\in\P_{\tau}(2n)$ and $m\in\N$,  we
define its $m$-th iteration $\ga^m:[0,m\tau]\to\Sp(2n)$ by
\be \ga^m(t) = \ga(t-j\tau)\ga(\tau)^j, \qquad
  {\rm for}\quad j\tau\leq t\leq (j+1)\tau,\;j=0,1,\ldots,m-1.
     \lb{3.6}\ee
We still denote the extended path on $[0,+\infty)$ by $\ga$.

{\bf Definition 3.2.} (cf. \cite{Lon3}, \cite{Lon5}) {\it For any $\ga\in\P_{\tau}(2n)$,
we define
\be (i(\ga,m), \nu(\ga,m)) = (i_1(\ga^m), \nu_1(\ga^m)), \qquad \forall m\in\N.
   \lb{3.7}\ee
The mean index $\hat{i}(\ga,m)$ per $m\tau$ for $m\in\N$ is defined by
\be \hat{i}(\ga,m) = \lim_{k\to +\infty}\frac{i(\ga,mk)}{k}. \lb{3.8}\ee
For any $M\in\Sp(2n)$ and $\om\in\U$, the {\it splitting numbers} $S_M^{\pm}(\om)$
of $M$ at $\om$ are defined by
\be S_M^{\pm}(\om)
     = \lim_{\ep\to 0^+}i_{\om\exp(\pm\sqrt{-1}\ep)}(\ga) - i_{\om}(\ga),
   \lb{3.9}\ee
for any path $\ga\in\P_{\tau}(2n)$ satisfying $\ga(\tau)=M$.}

For a given path $\gamma\in {\cal P}_{\tau}(2n)$ we consider to deform
it to a new path $\eta$ in ${\cal P}_{\tau}(2n)$ so that
\begin{equation}
i_1(\gamma^m)=i_1(\eta^m),\quad \nu_1(\gamma^m)=\nu_1(\eta^m), \quad
         \forall m\in {\bf N}, \label{3.10}
\end{equation}
and that $(i_1(\eta^m),\nu_1(\eta^m))$ is easy enough to compute. This
leads to finding homotopies $\delta:[0,1]\times[0,\tau]\to {\rm Sp}(2n)$
starting from $\gamma$ in ${\cal P}_{\tau}(2n)$ and keeping the end
points of the homotopy always stay in a certain suitably chosen maximal
subset of ${\rm Sp}(2n)$ so that (\ref{3.10}) always holds. In fact,  this
set was first discovered in \cite{Lon3} as the path connected component
$\Omega^0(M)$ containing $M=\gamma(\tau)$ of the set
\begin{eqnarray}
  \Omega(M)=\{N\in{\rm Sp}(2n)\,&|&\,\sigma(N)\cap{\bf U}=\sigma(M)\cap{\bf U}\;
{\rm and}\;  \nonumber\\
 &&\qquad \nu_{\lambda}(N)=\nu_{\lambda}(M)\;\forall\,
\lambda\in\sigma(M)\cap{\bf U}\}. \label{3.11}
\end{eqnarray}
Here $\Omega^0(M)$ is called the {\it homotopy component} of $M$ in
${\rm Sp}(2n)$.

In \cite{Lon3}-\cite{Lon5}, the following symplectic matrices were introduced
as {\it basic normal forms}:
\begin{eqnarray}
D(\lambda)=\left(\matrix{\lm & 0\cr
         0  & \lm^{-1}\cr}\right), &\quad& \lm=\pm 2,\lb{3.12}\\
N_1(\lm,b) = \left(\matrix{\lm & b\cr
         0  & \lm\cr}\right), &\quad& \lm=\pm 1, b=\pm1, 0, \lb{3.13}\\
R(\th)=\left(\matrix{\cos\th & -\sin\th\cr
        \sin\th  & \cos\th\cr}\right), &\quad& \th\in (0,\pi)\cup(\pi,2\pi),
                     \lb{3.14}\\
N_2(\om,b)= \left(\matrix{R(\th) & b\cr
              0 & R(\th)\cr}\right), &\quad& \th\in (0,\pi)\cup(\pi,2\pi),
                     \lb{3.15}\end{eqnarray}
where $b=\left(\matrix{b_1 & b_2\cr
               b_3 & b_4\cr}\right)$ with  $b_i\in\R$ and  $b_2\not=b_3$.

Splitting numbers possess the following properties:

{\bf Lemma 3.3.} (cf. \cite{Lon3} and Lemma 9.1.5 of \cite{Lon5}) {\it Splitting
numbers $S_M^{\pm}(\om)$ are well defined, i.e., they are independent of the choice
of the path $\ga\in\P_\tau(2n)$ satisfying $\ga(\tau)=M$ appeared in (\ref{3.9}).
For $\om\in\U$ and $M\in\Sp(2n)$, splitting numbers $S_N^{\pm}(\om)$ are constant
for all $N\in\Om^0(M)$. }

{\bf Lemma 3.4.} (cf. \cite{Lon3}, Lemma 9.1.5 and List 9.1.12 of \cite{Lon5})
{\it For $M\in\Sp(2n)$ and $\om\in\U$, there hold
\begin{eqnarray}
S_M^{\pm}(\om) &=& 0, \qquad {\it if}\;\;\om\not\in\sg(M).  \lb{3.16}\\
S_{N_1(1,a)}^+(1) &=& \left\{\matrix{1, &\quad {\rm if}\;\; a\ge 0, \cr
0, &\quad {\rm if}\;\; a< 0. \cr}\right. \lb{3.17}\eea

For any $M_i\in\Sp(2n_i)$ with $i=0$ and $1$, there holds }
\be S^{\pm}_{M_0\dm M_1}(\om) = S^{\pm}_{M_0}(\om) + S^{\pm}_{M_1}(\om),
    \qquad \forall\;\om\in\U. \lb{3.18}\ee

We have the following

{\bf Theorem 3.5.} (cf. \cite{Lon4} and Theorem 1.8.10 of \cite{Lon5}) {\it For
any $M\in\Sp(2n)$, there is a path $f:[0,1]\to\Om^0(M)$ such that $f(0)=M$ and
\be f(1) = M_1\dm\cdots\dm M_k,  \lb{3.19}\ee
where each $M_i$ is a basic normal form listed in (\ref{3.12})-(\ref{3.15})
for $1\leq i\leq k$.}

Let $\Sigma\in\H(2n)$. Using notations in \S1,
for any $(\tau,y)\in\J(\Sigma,\alpha)$ and $m\in\N$, we define
its $m$-th iteration $y^m:\R/(m\tau\Z)\to\R^{2n}$ by
\be y^m(t) = y(t-j\tau), \qquad {\rm for}\quad j\tau\leq t\leq (j+1)\tau,
       \quad j=0,1,2,\ldots, m-1. \lb{3.20}\ee
Note that this coincide with that in \S2.
We still denote by $y$ its extension to $[0,+\infty)$.

We define via Definition 3.2 the following
\bea  S^+(y) &=& S_{\ga_y(\tau)}^+(1),  \lb{3.21}\\
  (i(y,m), \nu(y,m)) &=& (i(\ga_y,m), \nu(\ga_y,m)),  \lb{3.22}\\
   \hat{i}(y,m) &=& \hat{i}(\ga_y,m),  \lb{3.23}\eea
for all $m\in\N$, where $\ga_y$ is the associated symplectic path of $(\tau,y)$.
Then we have the following.

{\bf Theorem 3.6.} (cf. Lemma 1.1 of \cite{LoZ1}, Theorem 15.1.1 of \cite{Lon5}) {\it Suppose
$(\tau,y)\in \J(\Sigma, \alpha)$. Then we have
\be i(y^m)\equiv i(m\tau ,y)=i(y, m)-n,\quad \nu(y^m)\equiv\nu(m\tau, y)=\nu(y, m),
       \qquad \forall m\in\N, \lb{3.24}\ee
where $i(y^m)$ and $\nu(y^m)$ are the index and nullity
defined in \S2. In particular, (\ref{2.7}) and (\ref{3.8})
coincide, thus we simply denote them by $\hat i(y)$.}

\setcounter{equation}{0}
\section{Proofs of the main theorems }

In the rest of this article, we fix a $\Sg\in\mathcal{SH}(2n)$ and assume the following
condition on $\Sg$:

\noindent (F) {\bf There exist only finitely many geometrically distinct
closed characteristics \\$\quad \{(\tau_j, y_j)\}_{1\le j\le k}$ on $\Sigma$. }

We denote by $\ga_j\equiv \gamma_{y_j}$ the associated
symplectic path of $(\tau_j,\,y_j)$ on $\Sg$ for $1\le j\le k$. Then by
Lemma 1.3 of \cite{LoZ1} or Lemma 15.2.4 of \cite{Lon5}, there exist $P_j\in \Sp(2n)$ and $M_j\in \Sp(2n-2)$ such
that
\be \ga_j(\tau_j)=P_j^{-1}(N_1(1,\,1)\dm M_j)P_j, \quad\forall\; 1\le j\le k,
   \lb{4.1}\ee
here we use notations in \S3.

Firstly we have the following property, cf., Lemma 4.2 of \cite{LLZ}.

{\bf Lemma 4.1.}  {\it Suppose $(\tau,\, y)\in\mathcal{J}(\Sigma, \,\alpha)$,
then $(\tau,\, -y)\in\mathcal{J}(\Sigma, \,\alpha)$
and either $\mathcal{O}(y)=\mathcal{O}(-y)$ or
$\mathcal{O}(y)\cap\mathcal{O}(-y)=\emptyset$,
where $\mathcal{O}(\pm y)=\{\pm y(t)|\, t\in\R\}$.
Moreover, if $\mathcal{O}(y)\cap\mathcal{O}(-y)\neq\emptyset$,
then we have}
$$ y(t)=-y\left(t+\frac{\tau}{2}\right),\qquad\forall t\in\R.$$

In the following we call a closed characteristic $(\tau,\, y)$
on $\Sg\in\mathcal{SH}(2n)$ {\it symmetric}
if $\mathcal{O}(y)\cap\mathcal{O}(-y)\neq\emptyset$,
{\it non-symmetric} if $\mathcal{O}(y)\cap\mathcal{O}(-y)=\emptyset$.
Thus if $(\tau,\, y)$ is  non-symmetric, then
$(\tau,\, y)$ and $(\tau,\, -y)$ are geometrically distinct;
if $(\tau,\, y)$ is symmetric, then
$(\tau,\, y)$ and $(\tau,\, -y)$ are geometrically the same.

We have the following property, cf., Lemma 15.6.4 of \cite{Lon5}.

{\bf Lemma 4.2.}  {\it Suppose $(\tau,\, y)\in\mathcal{J}(\Sigma, \,\alpha)$
is a symmetric closed characteristic on $\Sg\in\mathcal{SH}(2n)$. Then we have}
\be    i(y,\, 1)+2S^+(y)-\nu(y,\,1)\ge n. \lb{4.2}\ee

Now we can give:

{\bf Proof of Theorem 1.1.} Since $\hat i(y_j)>2$ for $1\le j\le k$, we can use
the common index jump theorem (Theorems 4.3 and 4.4 of
\cite{LoZ1}, Theorems 11.2.1 and 11.2.2 of \cite{Lon5}) to obtain infinitely many
$(T, m_1,\ldots,m_k)\in\N^{k+1}$ such that the following hold:
\bea
\nu(y_j,\, 2m_j-1) &=&\nu(y_j,\, 1), \lb{4.3}\\
i(y_j,\, 2m_j) &\ge& 2T-\frac{e(\gamma_j(\tau_j))}{2}\ge 2T-n, \lb{4.4}\\
i(y_j,\, 2m_j)+\nu(y_j,\, 2m_j) &\le& 2T+\frac{e(\gamma_j(\tau_j))}{2}-1\le 2T+n-1, \lb{4.5}\\
i(y_j,\, 2m_j+1) &=& 2T+i(y_j,\,1). \lb{4.6}\\
i(y_j,\, 2m_j-1)+\nu(y_j,\, 2m_j-1)
 &=& 2T-(i(y_j,\,1)+2S^+(y_j)-\nu(y_j, 1)). \lb{4.7}
\eea
Note that (\ref{4.5}) holds by Theorem 4.4 of \cite{LoZ1}, other parts
follows by Theorem 4.3 of \cite{LoZ1}.
More precisely, by Theorem 4.1 of
\cite{LoZ1} (in (11.1.10) in Theorem 11.1.1 of \cite{Lon5}, with $D_j=\hat i(y_j)$,
we have
\bea m_j=\left(\left[\frac{T}{M\hat i(y_j)}\right]+\chi_j\right)M,\quad 1\le j\le k,\lb{4.8}\eea
where $\chi_j=0$ or $1$ for $1\le j\le k$ and $M\in\N$ is fixed such that $\frac{M\theta}{\pi}\in\Z$,
whenever $e^{\sqrt{-1}\theta}\in\sigma(\ga_j(\tau_j))$ and $\frac{\theta}{\pi}\in\Q$
for some $1\le j\le k$. Moreover, we have the following
\be\frac{T}{M\hat i(y_j)}\in\N\quad{\rm and}\quad\chi_j=0, \qquad{\rm if}\quad \hat i(y_j)\in\Q,\lb{4.9}
\ee
which can be seen from the proof of  Theorem 4.1 of
\cite{LoZ1}, cf. the proof of Theorem 5.3 of \cite{LoZ1}.

By Corollary 1.2 of \cite{LoZ1}, we have
\be i(y_j,\,1)\ge n,\qquad 1\le j\le k.\lb{4.10}\ee
 Note that
$e(\gamma_j(\tau_j))\le2n$ for $1\le j\le k$. Hence Theorem 2.3
of \cite{LoZ1} yields
\bea i(y_j,\, m)+\nu(y_j,\, m)
&\le&  i(y_j, m+1)-i(y_j, 1)+\frac{e(\gamma_j(\tau_j))}{2}-1\nn\\
&\le& i(y_j, m+1)-1. \quad \forall m\in\N,\;1\le j\le k.\lb{4.11}\eea
Specially, we have
$$ i(y_j,\, m)<i(y_j,\, m+1),\qquad \forall m\in\N,\;1\le j\le k. $$

By Theorem 3.5 we have
\be  \gamma_j(\tau_j)\approx N_1(1,1)^{\diamond p_{j, -}}\diamond I_2^{\diamond p_{j, 0}}\diamond N_1(1,-1)^{\diamond  p_{j, +}}\diamond  G_j,
\qquad1\le j\le k\lb{4.12}\ee
for some nonnegative integers $p_{j, -}$, $p_{j, 0}$, $p_{j, +}$, and some symplectic
matrix $G_j$ satisfying $1\not\in \sigma(G_j)$.
By (\ref{4.1}), (\ref{4.12}) and Lemma 3.4 we obtain
\bea 2S^+(y_j)  = 2(p_{j, -} +p_{j, 0}) \ge 2,\quad1\le j\le k. \lb{4.13}
\eea

By (\ref{4.3}), (\ref{4.7}) and (\ref{4.10}) we have
\be i(y_j,\, 2m_j-1)
 = 2T-(i(y_j,\,1)+2S^+(y_j)\le 2T-n-2. \lb{4.14}
\ee
Now by (\ref{4.4})-(\ref{4.7}), (\ref{4.10}), (\ref{4.11}), (\ref{4.14}) and Theorem 3.6 we have
\bea
i(y_j^{2m_j}) &\ge& 2T-2n, \lb{4.15}\\
i(y_j^{2m_j})+\nu(y_j^{2m_j})-1 &\le& 2T-2, \lb{4.16}\\
i(y_j^{2m_j+m}) &\ge& 2T,\qquad \forall m\ge 1. \lb{4.17}\\
i(y_j^{2m_j-m})+\nu(y_j^{2m_j-m})-1
 &\le& 2T-2n-4,\qquad\forall m\ge 2. \lb{4.18}\\
i(y_j^{2m_j-1})+\nu(y_j^{2m_j-1})-1
 &=& 2T-(i(y_j,\,1)+2S^+(y_j)-\nu(y_j, 1))-n-1. \lb{4.19}
\eea
By Proposition 2.4, For every $1\le i\le n$, there exists a point
$u\in L_0^{\frac{\alpha}{\alpha-1}}(S^1,\R^{2n})$ such that
\bea
 \Phi^\prime(u)=0,\quad \Phi(u)=c_{T-i+1}, \quad
 C_{S^1,\; 2(T-i)}(\Phi, \;S^1\cdot u)\neq 0. \lb{4.20}\eea
Let $y_{\rho(i)}^{\lambda(i)}$ for $1\le i\le n$
be a set of closed characteristics satisfying (\ref{4.20}),
where $\rho: \{1,\ldots, n\}\rightarrow\{1,\ldots,k\}$ and
$\lambda: \{1,\ldots, n\}\rightarrow\N$ are integer valued functions.
Note that by Condition (F) and the infiniteness of the tuples
$(T, m_1,\ldots,m_k)\in\N^{k+1}$,
we can assume
\be y_{\rho(i)}\in\mathcal{V}_\infty(\Sg, \alpha),\qquad 1\le i\le n.\lb{4.21}\ee
By (\ref{4.17}), (\ref{4.18}), (\ref{4.20}) and Proposition 2.2, we have
\be \lambda(i)\in\{2m_{\rho(i)}-1,\;2m_{\rho(i)}\},\lb{4.22}\ee
for each $1\le i\le n$.

{\bf Claim 1.} {\it If $y_{\rho(i)}$ is symmetric, then $\lambda(i)=2m_{\rho(i)}$. }

In fact, by Lemma 4.2 and (\ref{4.19}), we have
\be i(y_{\rho(i)}^{2m_{\rho(i)}-1})+\nu(y_{\rho(i)}^{2m_{\rho(i)}-1})-1
 \le2T-2n-1. \lb{4.23}
\ee
Thus Claim 1 holds by (\ref{4.20}) and Proposition 2.2.

{\bf Claim 2.} {\it If $\lambda(i)=2m_{\rho(i)}-1$, then
$y_{\rho(i)}$ is non-symmetric and non-hyperbolic. }

The first statement follows directly from Claim 1. We prove the latter.
In fact, suppose $y_{\rho(i)}$ is hyperbolic.
Thus by (\ref{4.1}), (\ref{4.10}), Lemma 3.4 and Theorem 3.5 we have
\be    i(y_{\rho(i)},\, 1)+2S^+(y_{\rho(i)})-\nu(y_{\rho(i)},\,1)\ge n+1. \lb{4.24}\ee
Hence we  have
\be i(y_{\rho(i)}^{2m_{\rho(i)}-1})+\nu(y_{\rho(i)}^{2m_{\rho(i)}-1})-1
 \le2T-2n-2. \lb{4.25}
\ee
This contradict to (\ref{4.20}) and Proposition 2.2.
Thus Claim 2 holds.

{\bf Claim 3.} {\it If $\lambda(i_1)=2m_{\rho(i_1)}$ and
$\lambda(i_2)=2m_{\rho(i_2)}$, then either $\hat i(y_{\rho(i_1)})\in\R\setminus\Q$
or $\hat i(y_{\rho(i_2)})\in\R\setminus\Q$. }

Suppose the contrary, i.e., both $\hat i(y_{\rho(i_1)})\in\Q$ and $\hat i(y_{\rho(i_2)})\in\Q$.
Then by (\ref{4.8}) and (\ref{4.9}) we have
\bea  2m_{\rho(i_1)}\hat i(y_{\rho(i_1)})&=&2\left(\left[\frac{T}{M\hat i(y_{\rho(i_1)})}\right]+\chi_{\rho(i_1)}\right)M\hat i(y_{\rho(i_1)})
\nn\\&=&2\left(\frac{T}{M\hat i(y_{\rho(i_1)})}\right)M\hat i(y_{\rho(i_1)})=2T
=2\left(\frac{T}{M\hat i(y_{\rho(i_2)})}\right)M\hat i(y_{\rho(i_2)})
\nn\\&=&2\left(\left[\frac{T}{M\hat i(y_{\rho(i_2)})}\right]+\chi_{\rho(i_2)}\right)M\hat i(y_{\rho(i_2)})
=2m_{\rho(i_2)}\hat i(y_{\rho(i_2)}).
\lb{4.26}\eea
On the other hand, by (\ref{4.20}) we have
\be \Phi(y_{\rho(i_1)}^{2m_{\rho(i_1)}})=c_{T-i_1+1}\neq
c_{T-i_2+1}=\Phi(y_{\rho(i_2)}^{2m_{\rho(i_2)}}).\lb{4.27}\ee
By (\ref{4.21}) and Theorem 2.6 we have
\be \frac{\hat i(y_{\rho(i_1)})}{A(y_{\rho(i_1)})}=\frac{1}{\gamma(\Sigma)}=\frac{\hat i(y_{\rho(i_2)})}{A(y_{\rho(i_2)})}.\lb{4.28}
\ee
Note that we have the relations
\be \hat i(y^m)=m\hat i(y),\quad A(y^m)=mA(y),\quad
\Phi(y)=-\left(1-\frac{\alpha}{2}\right)\left(\frac{2}{\alpha}A(y)\right)^{\frac{\alpha}{\alpha-2}},\quad\forall m\in\N,\lb{4.29}\ee
for any closed characteristic $y$ on $\Sg$ by (\ref{2.7}), (\ref{2.12})
and (45) in P.221 of \cite{Eke3}.

Hence we have
\bea  2m_{\rho(i_1)}\hat i(y_{\rho(i_1)})
&=&\gamma(\Sigma)^{-1}\cdot 2m_{\rho(i_1)}A(y_{\rho(i_1)})
=\gamma(\Sigma)^{-1}\cdot A(y_{\rho(i_1)}^{2m_{\rho(i_1)}})\nn\\
&=&2(\gamma(\Sigma)C_\alpha)^{-1}(-\Phi(y_{\rho(i_1)}^{2m_{\rho(i_1)}})^{\frac{\alpha-2}{\alpha}}
=2(\gamma(\Sigma)C_\alpha)^{-1}(-c_{T-i_1+1})^{\frac{\alpha-2}{\alpha}}\nn\\
&\neq&2(\gamma(\Sigma)C_\alpha)^{-1}(-c_{T-i_2+1})^{\frac{\alpha-2}{\alpha}}
=2(\gamma(\Sigma)C_\alpha)^{-1}(-\Phi(y_{\rho(i_2)}^{2m_{\rho(i_2)}})^{\frac{\alpha-2}{\alpha}}\nn\\
&=&\gamma(\Sigma)^{-1}\cdot A(y_{\rho(i_2)}^{2m_{\rho(i_2)}})
=\gamma(\Sigma)^{-1}\cdot 2m_{\rho(i_2)}A(y_{\rho(i_2)})\nn\\
&=&2m_{\rho(i_2)}\hat i(y_{\rho(i_2)}),
\lb{4.30}\eea
where $C_\alpha$ is the constant given in Theorem 2.6.
This contradict to (\ref{4.26}).
Hence Claim 3 holds.

Now we prove Theorem 1.1 as follows.

For each $j\in \im\rho$, we have $^\#\{\rho^{-1}(j)\}\in\{1,\,2\}$
 by (\ref{4.22}). Then we have the following three cases.

{\bf Case 1.} { \it We have $j\in\Theta_1\equiv\{^\#\{\rho^{-1}(l)\}=2\}$.
}

In this case $y_j$ is non-symmetric and non-hyperbolic  by (\ref{4.22})
and Claim 2. Hence $y_j$ and $-y_j$ are geometrically distinct
by Lemma 4.1. Thus we obtain two non-hyperbolic closed characteristics
in $\td{\J}(\Sigma)$ for each $j$ such that $^\#\{\rho^{-1}(j)\}=2$.
Hence we have
$2^\#\Theta_1$ non-hyperbolic closed characteristics
in $\td{\J}(\Sigma)$ in this case.

{\bf Case 2.} { \it We have $j\in\Theta_2\equiv\{^\#\{\rho^{-1}(l)\}=1\; {\rm and}
\;\lambda(\rho^{-1}(l))=2m_l-1\}$.}

In this case $\rho|_{ \rho^{-1}(\Theta_2)}: \rho^{-1}(\Theta_2)\rightarrow\Theta_2$
is a bijection. By Claim 2,
$y_j$ is non-symmetric and non-hyperbolic for $j\in\Theta_2$.
Thus we obtain two non-hyperbolic closed characteristics
for each $j\in\Theta_2$ as in Case 1. Hence we have
$2^\#\Theta_2$ non-hyperbolic closed characteristics
in $\td{\J}(\Sigma)$ in this case.

{\bf Case 3.} { \it We have $j\in\Theta_3\equiv\{^\#\{\rho^{-1}(l)\}=1\; {\rm and}
\;\lambda(\rho^{-1}(l))=2m_l\}$.}

Note that in this case $\rho|_{ \rho^{-1}(\Theta_3)}: \rho^{-1}(\Theta_3)\rightarrow\Theta_3$
is a bijection. By Claim 3, there exists at most one $j\in\Theta_3$
such that $\hat i(y_j)\in\Q$. Thus there are at least $^\#\Theta_3-1$
non-hyperbolic closed characteristics in $\td{\J}(\Sigma)$ in this case.

Clearly $\Theta_1, \,\Theta_2, \,\Theta_3$
are pairwise disjoint and
$\{1,\ldots, n\}=\rho^{-1}(\Theta_1)\cup\rho^{-1}(\Theta_2)\cup\rho^{-1}(\Theta_3)$.
Thus we have $2^\#\Theta_1+^\#\Theta_2+^\#\Theta_3=n$ since
$\rho|_{ \rho^{-1}(\Theta_1)}: \rho^{-1}(\Theta_1)\rightarrow\Theta_1$
is a two to one map and $\rho|_{ \rho^{-1}(\Theta_i)}: \rho^{-1}(\Theta_i)\rightarrow\Theta_i$
are bijections for $i=2,\,3$. By Claims 1-3, the number of
non-hyperbolic closed characteristics in $\td{\J}(\Sigma)$
is at least
$$2^\#\Theta_1+2^\#\Theta_2+^\#\Theta_3-1\ge 2^\#\Theta_1+^\#\Theta_2+^\#\Theta_3-1=n-1.$$
The proof of Theorem 1.1 is complete.
\hfill\hb

{\bf Proof of Theorem 1.2.}
As Definition 1.1 of \cite{LoZ1}, for $\alpha\in(1,2)$, we define a map
$\varrho_n\colon\H(2n)\to\N\cup\{ +\infty\}$
\be \varrho_n(\Sg)
= \left\{\matrix{+\infty, & {\rm if\;\;}^\#\mathcal{V}(\Sigma,\alpha)=+\infty, \cr
\min\left\{[\frac{i(x,1) + 2S^+(x) - \nu(x,1)+n}{2}]\,
\left|\frac{}{}\right.\,(\tau,x)\in\mathcal{V}_\infty(\Sigma, \alpha)\right\},
 & {\rm if\;\;} ^\#\mathcal{V}(\Sigma, \alpha)<+\infty, \cr}\right.  \lb{4.31}\ee
where $\mathcal{V}(\Sigma,\alpha)$ and $\mathcal{V}_\infty(\Sigma,\alpha)$ are
variationally visible and infinite variationally visible sets respectively given
in Definition 2.5.

By Theorem 1.4 of \cite{LoZ1}, if $^\#\td{\J}(\Sigma)\le 2\varrho_n(\Sg)-2$,
then there exist at least two elliptic closed characteristics
in $\td{\J}(\Sigma)$. By Theorem 1.1 of \cite{LoZ1} we have
$\varrho_n(\Sg)\ge [\frac{n}{2}]+1$. Thus Theorem 1.2 holds
when $n$ is even.

In the following, we prove Theorem 1.2 for $n$ being odd.
We have the following two cases.

{\bf Case 1.} { \it All the closed characteristics on $\Sg$ are symmetric.}

In this case, by Lemma 4.2 and (\ref{4.31}), we have $\varrho_n(\Sg)\ge n$.
Thus Theorem 1.2 holds by Theorem 1.4 of \cite{LoZ1}.

{\bf Case 2.} { \it At least one closed characteristic on $\Sg$ is non-symmetric.}

We may assume without loss of generality that $(\tau_2,\, y_2)=(\tau_1,\, -y_1)$.
Since $\Sg=-\Sg$, we have $H_\alpha(x)=H_\alpha(-x)$.
Then it is easy to see that $(\tau_1, y_1)$ and $(\tau_1, -y_1)$
have the same properties
\bea&&(i(y_1^m),\,\nu(y_1^m))=(i((-y_1)^m),\,\nu((-y_1)^m)),\quad
\Phi(y_1^m)=\Phi((-y_1)^m),\qquad\forall m\in\N\lb{4.32}\\
&&C_{S^1,\; q}(\Phi, \;S^1\cdot u_1^m)\cong C_{S^1,\; q}(\Phi, \;S^1\cdot (-u_1)^m),\qquad \forall m\in\N,\;\forall q\in\Z,\lb{4.33}
\eea
where we denote by $(\pm u)^m$ the critical point
of $\Phi$ corresponding to $(\pm y)^m$. In fact,
we have a nature $\Z_2$-action on
$L_0^{\frac{\alpha}{\alpha-1}}(S^1,\R^{2n})$ defined by
$u\mapsto-u$ and the functional $\Phi$ defined in (\ref{2.3})
is $\Z_2$-invariant. Thus (\ref{4.32}) and (\ref{4.33}) hold.

Now we consider the set of closed characteristics:
$\Delta\equiv\{(\tau_1,\, y_1),\;(\tau_3,\, y_3),\ldots,(\tau_n,\,y_n\}$,
i.e., we remove $(\tau_2,\, y_2)$ from the set
$\{(\tau_j, y_j)\}_{1\le j\le n}$.
Then one can use the proof of Theorem 1.4 of \cite{LoZ1}
to obtain Theorem 1.2. In fact, we have $^\#\Delta=n-1=2[\frac{n}{2}]$.
While the proof of Theorem 1.4 of \cite{LoZ1} depend only on the
index iteration theory, hence it remains valid
if we replace the set $\{(\tau_j, y_j)\}_{1\le j\le n}$
there by $\Delta$. The proof of Theorem 1.2 is complete.\hfill\hb

\noindent {\bf Acknowledgements.} I would like to sincerely thank my
advisor, Professor Yiming Long, for introducing me to Hamiltonian
dynamics and for his valuable helps to me in all ways.
I would like to say that how enjoyable it is to work with him.

\bibliographystyle{abbrv}

\medskip

\end{document}